\documentclass[10pt]{amsart}
\usepackage{amsmath}
  \usepackage{paralist}
  \usepackage{graphics} 

\usepackage{latexsym,amssymb}
\usepackage{hyperref}
\usepackage{amsmath,amsthm}
\usepackage{amsfonts}
\usepackage[american]{babel}
\usepackage{mathrsfs}
\usepackage{psfrag}

\newtheorem{theorem}{Theorem}[section]

\newtheorem{lemma}[theorem]{Lemma}

\theoremstyle{definition}

\newtheorem{remark}{Remark}

\newcommand{\R}{\mathbb{R}}

\newcommand{\C}{\mathbb{C}}

\title[Symmetries for the Green's function]
      {Symmetries in an overdetermined problem \\ for the Green's function}

\author[Virginia Agostiniani and Rolando Magnanini]{}


 \email{vagostin@sissa.it}
 \email{magnanin@math.unifi.it}

\begin{document}
\maketitle

\centerline{\scshape Virginia Agostiniani}
\medskip
{\footnotesize
 \centerline{SISSA}
 \centerline{via Beirut 2-4, 34151 Trieste, ITALIA}
} 

\medskip

\centerline{\scshape Rolando Magnanini}
\medskip
{\footnotesize
 \centerline{Dipartimento di Matematica ``U. Dini''}
 \centerline{Universit\`a degli Studi di Firenze}
 \centerline{viale Morgagni 67/A, 50134 Firenze, ITALIA}
}

\bigskip


\begin{abstract}
We consider in the plane the problem of reconstructing a domain from the
normal derivative of its Green's function with pole at a fixed point in the domain. 
By means of the theory of conformal mappings, we obtain existence, uniqueness, (non-spherical) symmetry results,
and a formula relating the curvature of the boundary of the domain to the normal derivative of its Green's function.
\end{abstract}

\section{Introduction}

Overdetermined boundary value problems in partial differential equations have connections
to various fields in mathematics; they emerge in the study of isoperimetric inequalities,
optimal design and ill-posed and free boundary problems, to name a few. In many such
problems one's interest is focused on a specific feature: the shape of the domain considered;
mainly, its (spherical) symmetry, as in Serrin's landmark paper \cite{Se} and its many offsprings
(see \cite{We}, \cite{Al}, \cite{Fr}, \cite{Le}, \cite{Pa}, and the references therein).

With the present paper, we want to start a more detailed analysis of overdetermined problems in the plane,
by exploiting the full power of the theory of analytic functions.
As a case study, we shall analyse what appears to be the simplest situation:
in a planar bounded domain $\Omega$ with boundary $\partial\Omega$
of class $C^{1,\alpha}$, we shall consider the problem
\begin{eqnarray}
                   -\Delta U\!\!\!&=&\!\!\!\delta_{\zeta_c} \hspace{1.5cm}  \mbox{in } \Omega, \label{pb1} \\
                   U\!\!\!&=&\!\!\!0 \hspace{1.5cm} \mbox{on } \partial\Omega, \label{pb2} \\
                   \frac{\partial U}{\partial\nu}\!\!\!&=&\!\!\!\varphi
                   \hspace{1.5cm}\mbox{on } \partial\Omega. \label{pb3}
\end{eqnarray}
where $\nu$ is the {\it interior} normal direction to $\partial\Omega$,
$\delta_{\zeta_c}$ is the Dirac delta centered at a given point $\zeta_c\in\Omega$
and $\varphi:\partial\Omega\rightarrow\R$ is a positive given function of arclength,
measured counterclockwise from a reference point on $\partial\Omega$.

Problem (\ref{pb1})-(\ref{pb3}) can be interpreted as a free-boundary problem: find
a domain $\Omega$ whose Green's function $U$ with pole at $\zeta_c$ has gradient
with values on the boundary that fit those of the given function $\varphi$.
This formulation serve as a basis to model, for example, the Hele-Shaw flow, as
done in \cite{Gu} and \cite{Sa}.

By means of the Riemann Mapping Theorem, the solution of (\ref{pb1})-(\ref{pb2})
can be esplicitly written in terms of a conformal mapping $f$ from the unit disk $D$
to $\Omega$, which is uniquely determined if satisfies some suitable normalizing conditions.
Since it turns out that the normal derivative of $U$ on $\partial\Omega$ is proportional
to the modulus of the inverse of $f,$ 
then by (\ref{pb3}) and classical results on holomorphic functions, we can derive an explicit formula
for $f$ in terms of $\varphi$ (see section \S2 for details). With the help
of such a formula, we obtain the following results:
\begin{itemize}

\item[(i)] existence and uniqueness theorems for a domain $\Omega$
satisfying (\ref{pb1})-(\ref{pb3}) (Theorems \ref{tom} and \ref{esistenza}); 

\item[(ii)] symmetry results relating the invariance of $\varphi$ under certain groups
of transformations to that of $\Omega$ (Theorems \ref{rot} and \ref{simasse});
of course, when $\varphi$ in constant,
we obtain that $\Omega$ is a disk --- a well-known result (see \cite{Pa}, \cite{Le} \cite{Al});

\item[(iii)] a formula relating the interior normal derivative of the
Green's function to the curvature of $\partial\Omega.$

\end{itemize}

\section{Construction of a forward operator and its inverse}

In what follows, $D$ will always be the open unit disk in $\C$ centered at $0.$ 

Let us recall some basic facts of harmonic and complex analysis. We refer the reader to
\cite{go} and \cite{Ma} for more details. If $\Omega\subseteq\C$ is a simply connected domain 
bounded by a Jordan curve and
$\zeta_c\in\Omega$, then, by the Riemann Mapping Theorem, $\Omega$ is the image of an analytic function
$f:D\rightarrow\Omega$ which induces a homeomorphism between the closures
$\overline D$ and $\overline\Omega$, has non-zero derivative $f'$ in $D$ and
is such that $f(0)=\zeta_c.$
Moreover, if $\Omega$ is of class $C^{1,\alpha},$ $0<\alpha<1,$ that is its boundary 
$\partial\Omega$ is locally the graph of a function of class $C^{1,\alpha},$ then,
by Kellogg's theorem, we can infer that $f\in C^{1,\alpha}(\overline D)$ (see \cite{go}).

The following elementary lemma will be useful in the sequel.

\begin{lemma}\label{teo.}
Let $\Omega$ be a bounded simply connected domain in $\C$ and
$f:D\rightarrow\Omega$ be one-to-one and analytic with $f\in C^1(\overline D),$ $0<\alpha<1.$
Then there exists $\gamma\in\R$ such that
\begin{equation}\label{rappre}
f'(z)=e^{i\gamma}\exp\left\{\frac 1{2\pi}
\int_0^{2\pi}\frac{e^{it}+z}{e^{it}-z}\log|f'(e^{it})|dt\right\}
\end{equation}
for every $z\in D$.
\end{lemma}

\begin{proof}
The function
$$
f'(z)\exp\left\{-\frac 1{2\pi}
\int_0^{2\pi}\frac{e^{it}+z}{e^{it}-z}\log|f'(e^{it})|dt\right\},
\hspace{1cm}z\in D,
$$
is analytic, never zero in $D$ and has unitary modulus on $\partial D$; hence
it equals the number $e^{i\gamma}$ for some $\gamma\in\R.$
\end{proof}

With these premises, given two distinct numbers $\zeta_c$ and $\zeta_b\in\C,$ we consider 
\begin{center}
the set 
${\mathscr O}$ of all $C^{1,\alpha},$ $0<\alpha<1,$
simply connected
\end{center}
\begin{center} 
bounded domains such that
$\zeta_c\in\Omega$ and $\zeta_b\in\partial\Omega$.
\end{center}
We can put $\mathscr O$ in one-to-one
correspondence with 
\begin{center}
the class
$\mathscr F$ of all one-to-one analytic mappings 
\end{center}
\begin{center}
$f\in C^{1,\alpha}(\overline D)$
such that $f(0)=\zeta_c$ and $f(1)=\zeta_b$.
\end{center} 
In fact, the arbitrary parameter
$\gamma$ in (\ref{rappre}) can be determined by observing that 
\begin{equation}\label{detalfa}
\zeta_b-\zeta_c=\int_0^1f'(t)dt.
\end{equation}
\par
We now construct our forward operator $\mathcal{T}$ as the one that associates to each $\Omega$ in $\mathscr O$ 
the interior normal derivative $\frac{\partial U}{\partial\nu}$ --- as function of the arclength,
measured counterclockwise on $\partial\Omega,$ starting from $\zeta_b$ --- of the solution
of (\ref{pb1})-(\ref{pb2}). With our identification of $\mathscr O$ with $\mathscr F$ in mind, for
$f\in\mathscr F$, $\mathcal{T}(f)$ is a function of arclength $s\in[0,|\partial\Omega|]$ and it is defined
by the following remarks.

First, notice that, by Gauss-Green's formula, if $U$ satisfies (\ref{pb1})-(\ref{pb2}), then
\begin{equation*}
v(\zeta_c)=\int_{\partial\Omega}v(\zeta)\frac{\partial U}{\partial\nu}(\zeta)ds(\zeta)
\end{equation*}
for every function $v\in C^1(\overline\Omega)\cap C^2(\Omega)$ which is harmonic in $\Omega.$

Secondly, recall that any such function $v$ satisfies the well-known formula 
\begin{equation*}
v(\zeta)=\frac 1{2\pi}\int_{\partial\Omega}v(\zeta')
\frac{1-|f^{-1}(\zeta)|^2}{|f^{-1}(\zeta)-f^{-1}(\zeta')|}\frac{ds(\zeta')}{|f'(f^{-1}(\zeta'))|},
\hspace{1cm}\zeta\in\Omega,
\end{equation*}
if $\partial\Omega$ is rectifiable (see \cite{Ma}). By comparing the last two formulas (with $\zeta=\zeta_c=f(0)),$ we obtain that
\begin{equation*}
\frac{\partial U}{\partial\nu}(\zeta)=\frac 1{2\pi|f'(f^{-1}(\zeta))|},
\hspace{1.5cm}\zeta\in\partial\Omega.
\end{equation*}

Thirdly, since the arclength on $\partial\Omega$ is related to $f$ by the formula 
\begin{equation}
\label{arclength}
s(\theta)=\int_0^{\theta}|f'(e^{it})|dt,\hspace{1.5cm}\theta\in[0,2\pi],
\end{equation}
the values $\mathcal{T}(f)(s),$ $s\in[0,|\partial\Omega|],$ can be defined parametrically
by
\begin{equation}\label{param}
s = \int_0^{\theta}|f'(e^{it})|dt,\ \ \mathcal T(f)=\frac 1{2\pi|f'(e^{i\theta})|},
\ \ \theta\in[0,2\pi].
\end{equation}
It is clear that $\mathcal T(f)\in C^{0,\alpha}([0,|\partial\Omega|])$ 
and also that 
$$
\int_0^{|\partial\Omega|}\mathcal{T}(f)(s)ds=1,\hspace{1.5cm}\mathcal{T}(f)>0\mbox{ on }[0,|\partial\Omega|],
$$
for all $f\in\mathscr F.$

We shall now prove that $\mathcal T$ is injective by showing
that each $\varphi$ in the range of $\mathcal{T}$ determines only one $f\in\mathscr F.$  
In fact, for
$\varphi\in C^{0,\alpha}([0,|\partial\Omega|])$ in the 
range of $\mathcal{T}$, by formulas (\ref{param}) it turns out that
\begin{equation}\label{fi}
2\pi\varphi(s(\theta))s'(\theta)=1,\hspace{1.5cm}\theta\in[0,2\pi].
\end{equation}
This last formula, once
integrated between $0$ and $\theta$, gives
\begin{equation}\label{teta}
s(\theta)=\Phi^{-1}(\theta),\hspace{1.5cm}\theta\in[0,2\pi],
\end{equation}
where $\Phi^{-1}$ is the inverse of $\Phi:[0,|\partial\Omega|]\rightarrow[0,2\pi]$ defined by
\begin{equation}\label{teta2}
\Phi(s)=2\pi\int_0^s\varphi(\sigma)d\sigma,\hspace{1.5cm}s\in[0,|\partial\Omega|].
\end{equation}
By the same formulas (\ref{param}), we then obtain that 
\begin{equation}\label{f' g}
|f'(e^{i\theta})|=\frac 1{2\pi\varphi(\Phi^{-1}(\theta))},\hspace{1.5cm}\theta\in[0,2\pi],
\end{equation}
and hence (\ref{rappre}) gives
\begin{equation}\label{rappre2}
f'(z)=e^{i\gamma}\exp\left\{\frac 1{2\pi}
\int_0^{2\pi}\frac{e^{it}+z}{e^{it}-z}\log\frac 1{2\pi\varphi(\Phi^{-1}(t))}dt\right\},
\hspace{1.5cm}z\in D,
\end{equation}
where $\gamma$ is determined by (\ref{detalfa}).
Therefore, for any $\varphi$ in the range of $\mathcal{T}$, a unique $f\in\mathscr F$
such that $\mathcal{T}(f)=\varphi$  is determined by
$$
f(z)=\zeta_c+\int_0^1f'(tz)zdt,\hspace{1.5cm}z\in D,
$$
with $f'$ given by (\ref{rappre2}).

We collect these remarks in the following theorem.

\begin{theorem}\label{tom}
Given $\Omega\in\mathscr O$, let $\zeta_b$ be a reference point on $\partial\Omega$ from which the
arclength on $\partial\Omega$ is measured counterclockwise.  
\par
If $\varphi$ is the interior normal derivative of the Green's function on $\partial\Omega$ (as function of the
arclength), then a unique $f\in\mathscr F$ is determined such that
$\mathcal{T}(f)=\varphi$ and its derivative is given by
\begin{equation}\label{nicola}
f'(z)=e^{i\gamma}\exp\left\{\frac 1{2\pi}
\int_0^{2\pi}\frac{e^{it}+z}{e^{it}-z}\log\frac 1{2\pi\varphi(s(t))}dt\right\},
\ \ z\in D,
\end{equation}
where $s$ and $\Phi$ are defined by (\ref{teta}) and (\ref{teta2}), respectively.

Moreover, the constant $\gamma$ is determined by
\begin{equation}\label{condi}
e^{i\gamma}\int_0^1\exp\left\{\frac 1{2\pi}\int_0^{2\pi}\frac{e^{i\tau}+t}{e^{i\tau}-t}
\log\frac 1{2\pi\varphi(s(\tau))}d\tau\right\}dt=\zeta_b-\zeta_c.
\end{equation}
\end{theorem}

\vspace{.5cm}

Theorem \ref{tom} tells us that the operator $\mathcal T$ is injective. A discussion about
its surjectivity is beyond the aims of this paper.  
Far from being complete, we want here to suggest the
following criterion.

Referring to \cite{Du}, let us introduce the so called \emph{boundary rotation} of a function $f$
defined in $D$:
$$
\rho=\lim_{r\rightarrow1^-}
\int_0^{2\pi}\left|{\mbox Re}\left\{1+\frac{zf''(z)}{f'(z)}\right\}
\right|d\theta,\ \ z=re^{i\theta}\in D.
$$
We consider the class $\mathscr V$ of all normalized functions
$$
f(z)=z+a_2z^2+a_3z^3+...
$$
which are analytic, locally univalent and with $\rho<+\infty$.
The proof of the surjectivity of $\mathcal T$ relies 
on the problem of finding an analytic and univalent function
$f$ from the disk to $f(D)=\Omega$. Therefore, we have to look for sufficient conditions for univalence. 
The following theorem is based on a sufficient condition, due to Paatero, that
says that any function in 
the class $\mathscr V$ with $\rho\le4\pi$ is univalent (see \cite{Du}). 

\begin{theorem}\label{esistenza}
Let $\varphi\in C^1(\R)$ be $L$-periodic, strictly positive and satisfying the compatibility condition
$\int_0^L\varphi(s)ds=1.$ If, moreover, $\varphi$ satisfies the condition
$$
\max_{[0,L]}\left|\frac{\varphi'(s)}{\varphi^2(s)}\right|\le2\pi,
$$
then $\mathcal T$ is surjective.
\end{theorem} 

\begin{proof}
By Theorem \ref{tom}, we know that a function $f\in\mathscr F$ such that $\mathcal T(f)=\varphi$ must satisfy
(\ref{nicola}). Thus, we have to check Paatero's condition on (\ref{nicola}). From that expression
we deduce that
$$
\frac{f''(z)}{f'(z)}=\frac1{2\pi}\int_0^{2\pi}\frac{2e^{it}}{(e^{it}-z)^2}\log\frac1{2\pi\varphi(s(t))}dt,
$$ 
being $s$ defined as in (\ref{teta}) and (\ref{teta2}).
Now, by observing that
$$
\frac{d}{dt}\left(\frac{e^{it}+z}{e^{it}-z}\right)=\frac{-2ize^{it}}{(e^{it}-z)^2},
$$
we can integrate by parts and obtain that
$$
\frac{-izf''(z)}{f'(z)}=\frac1{2\pi}\int_0^{2\pi}\frac{e^{it}+z}{e^{it}-z}\frac{\varphi'(s(t))s'(t)}
{\varphi(s(t))}dt.
$$
By the maximum modulus principle, we can estimate, for $z\in D,$
\begin{eqnarray*}
\left|{\mbox Re}\left\{1+
\frac{zf''(z)}{f'(z)}
\right\}\right|        &\le& 1+\left|\frac{-izf''(z)}{f'(z)}\right| \\
                       &\le& 1+\max_{[0,2\pi]}\left|\frac{\varphi'(s(t))s'(t)}{\varphi(s(t))}\right|,                             
\end{eqnarray*}
and, from (\ref{fi}), we have that $\varphi'(s)s'/\varphi(s)=\varphi'(s)/ 2\pi\varphi^2(s).$
Therefore, we can estimate the boundary rotation of $f$ in the following way:
$$
\rho\le\int_0^{2\pi}\left(1+\max_{[0,2\pi]}
\left|\frac{\varphi'(s(t))}{2\pi\varphi^2(s(t))}\right|\right)d\theta=
2\pi\left(1+\max_{[0,L]}
\left|\frac{\varphi'(s)}{2\pi\varphi^2(s)}\right|\right).
$$
By our assumptions, it follows that $\rho\le4\pi$ and hence, from Paatero's
criterion for univalence, $f$ is a homeomorphism from the disk onto $f(D).$
\end{proof}

\section{Symmetries}

\begin{remark}
Theorem \ref{tom} allows us to rediscover a result already proved in \cite{Pa} and also in
\cite{Le} and \cite{Al}:
if $\varphi$ is constant, then $\Omega$ is a disk. More precisely,
given $\Omega\in\mathscr O$ with perimeter $L$, let
$\varphi$ be constantly equal to $C>0.$ 
From (\ref{nicola}), we obtain that
$$
f'(z)=e^{i\gamma}\exp\left\{\frac 1{2\pi}\log\frac 1{2\pi C}
\int_0^{2\pi}\frac{e^{it}+z}{e^{it}-z}dt\right\}=\frac{e^{i\gamma}}{2\pi C},
$$
since $\int_0^{2\pi}\frac{e^{it}+z}{e^{it}-z}dt=2\pi$. Therefore,
we get that
$$
f(z)=\zeta_c+\frac{e^{i\gamma}}{2\pi C}\, z,\hspace{1.5cm}z\in D,
$$
that is $\Omega$ is the disk centered at $\zeta_c$ with radius $\frac 1{2\pi C}.$
\end{remark}

Now we want to show how some other symmetry properties of $\Omega$ 
can be derived from some invariance properties of $\varphi$ 
and viceversa. 

In what follows, for $\Omega\in\mathscr O,$ let $L=|\partial\Omega|$ and let $\varphi$ denote 
the values of the interior normal
derivative on $\partial\Omega$ (as function of arclength) of the Green's function
of $\Omega.$ 

In the next theorem, we will identify $\varphi$ 
with its $L$-periodic extension to $\R$ and $\mathcal{R}_{\zeta,\beta}$ will denote the clockwise
rotation of an angle $\beta$
around a point $\zeta.$


\begin{theorem}\label{rot}
Given $\Omega\in\mathscr O,$
for every $n=2,3,...,$
\begin{center}
$\mathcal{R}_{\zeta_c,\frac{2\pi}n}(\Omega)=\Omega$
if and only if $\varphi$ is $\frac Ln$-periodic.
\end{center}
\end{theorem}

\begin{proof}
Let us fix $n$ and suppose $\varphi$ measured counterclockwise from
$\zeta_b\in\partial\Omega$. Let $f\in\mathscr F$ be the unique analytic function from $D$ to $\Omega$
such that $f(0)=\zeta_c$ and $f(1)=\zeta_b.$

(i) If $\Omega$ is invariant by rotations of angle $\frac{2\pi}n$ around $\zeta_c,$ 
then $f$ satisfies
$$
f(ze^{i\frac {2\pi}n})=\zeta_c+[f(z)-\zeta_c]e^{i\frac {2\pi}n},\hspace{1.5cm}z\in D.
$$
By differentiating this expression, we obtain $f'(ze^{i\frac {2\pi}n})=f'(z),$ from which
$$
s\left(\theta +\frac{2\pi}n\right)=\int_{0}^{\theta+\frac{2\pi}n}|f'(e^{it})|dt=
s(\theta)+\int_{\theta}^{\theta+\frac{2\pi}n}|f'(e^{it})|dt,
$$
and hence
\begin{equation}
\label{stheta}
s\left(\theta +\frac{2\pi}n\right)=
s(\theta)+s\left(\frac{2\pi}n\right),\ \ \theta\in\R.
\end{equation}
Now, being
$$
L=s(2\pi)=s\left(\frac{n-1}n2\pi\right)
+s\left(\frac{2\pi}n\right)=...=ns\left(\frac{2\pi}n\right),
$$
we have that $s\left(\theta+\frac{2\pi}n\right)=s(\theta)+\frac Ln.$ 
Thus, \eqref{stheta} and \eqref{fi}-\eqref{f' g} imply that
$$
\varphi\left(s(\theta)+\frac Ln\right)=\varphi\left(s\left(\theta+\frac{2\pi}n\right)\right)
=\frac 1{2\pi|f'(e^{i(\theta+\frac{2\pi}n)})|}=\frac 1{2\pi|f'(e^{i\theta})|}=
\varphi(s(\theta)),
$$
and hence, for every $s\in\R,$
$$
\varphi\left(s+\frac Ln\right)=\varphi(s).
$$
(ii) If now $\varphi$ is $\frac Ln$-periodic, from (\ref{teta2}) we write
\begin{equation}\label{shu}
\Phi\left(s+\frac Ln\right)=2\pi\int_0^{s+\frac Ln}\varphi(\sigma)d\sigma=
\Phi(s)+\Phi\left(\frac Ln\right).
\end{equation}
Since (\ref{teta}) holds, it follows that
$$
2\pi=\Phi(s(2\pi))=\Phi(L)=\Phi\left(\frac{n-1}nL\right)+
\Phi\left(\frac Ln\right)=...=n\Phi\left(\frac Ln\right),
$$
and hence
$$
\Phi\left(\frac Ln\right)=\frac{2\pi}n=\theta+\frac{2\pi}n-\theta=
\Phi\left(s\left(\theta+\frac{2\pi}n\right)\right)-\Phi(s(\theta)).
$$
From this and (\ref{shu}), we infer that
$$
\Phi\left(s\left(\theta+\frac{2\pi}n\right)\right)=
\Phi(s(\theta))+\Phi\left(\frac Ln\right)=\Phi\left(s(\theta)+\frac Ln\right),
$$
and, thanks to the invertibility of $\Phi$, we obtain
$$
s\left(\theta+\frac{2\pi}n\right)=s(\theta)+\frac Ln,\hspace{1.5cm}\theta\in\R.
$$
By this formula, (\ref{nicola}) and the periodicity of $\varphi$, it follows that
\begin{eqnarray*}
f'(z)&=& e^{i\gamma}\exp\left\{\frac 1{2\pi}
         \int_0^{2\pi}\frac{e^{it}+z}{e^{it}-z}\log\frac 1{2\pi\varphi(s(t))}dt\right\} \\
     &=& e^{i\gamma}\exp\left\{\frac 1{2\pi}
         \int_0^{2\pi}\frac{e^{it}+z}{e^{it}-z}\log\frac 1{2\pi\varphi\left(s\left(\frac{2\pi}n+t\right)-\frac Ln\right)}dt\right\} \\
     &=& e^{i\gamma}\exp\left\{\frac 1{2\pi}
         \int_0^{2\pi}\frac{e^{it}+z}{e^{it}-z}\log\frac 1{2\pi\varphi\left(s\left(\frac{2\pi}n+t\right)\right)}dt\right\}.
\end{eqnarray*}
By a change of variables, we thus get
\begin{eqnarray*}
f'(z)&=& e^{i\gamma}\exp\left\{\frac 1{2\pi}
         \int_{\frac{2\pi}n}^{2\pi+\frac{2\pi}n}
         \frac{e^{i(t-\frac{2\pi}n)}+z}{e^{i(t-\frac{2\pi}n)}-z}
         \log \frac 1{2\pi\varphi(s(t))}dt\right\} \\
     &=& e^{i\gamma}\exp\left\{\frac 1{2\pi}
         \int_0^{2\pi}
         \frac{e^{it}+ze^{i\frac{2\pi}n}}{e^{it}-ze^{i\frac{2\pi}n}}
         \log \frac 1{2\pi\varphi(s(t))}dt\right\} \\
     &=& f'(ze^{i\frac{2\pi}n}).
\end{eqnarray*}
Finally we find
$$
f(z)-\zeta_c=\int_0^1f'(tz)zdt=\int_0^1f'(tze^{i\frac{2\pi}n})zdt=[f(ze^{i\frac{2\pi}n})-\zeta_c]
e^{-i\frac{2\pi}n},
$$
and hence $\mathcal{R}_{\zeta_c,\frac{2\pi}n}\Omega=\Omega$.
\end{proof}

%
%

In what follows, $\mathcal{M}$ will denote mirror-reflection with
respect to a given axis. 


\begin{theorem}\label{simasse}
A domain $\Omega\in\mathscr O$ 
is symmetric with respect
to a generic axis if and only if
$\varphi(s)=\varphi(L-s)$ for all $s\in[0,L].$
\par
Here $\varphi$ is measured counterclockwise starting from an intersection point of the
axis with $\partial\Omega$.
\end{theorem}

\begin{proof} 
(i) Suppose $\Omega$ symmetric with rispect
to a given axis, that is $\mathcal{M}(\Omega)=\Omega.$
Short of rotations and translations, we can assume the symmetry axis to coincide with the real axis,
so that $\mathcal{M} z$ is the conjugate $\overline{z}$ of $z.$ 
\par
Let $f\in\mathscr F$ be the unique mapping from $D$ to $\Omega$
such that $f(0)=\zeta_c$ and $f(1)=\zeta_b,$ where $\zeta_c\in\Omega$ and
$\zeta_b\in\partial\Omega$ are some reference points on the symmetry axis.
We keep in mind that arclength on $\partial \Omega$ is measured counterclockwise from $\zeta_b.$
\par
It is clear that $\zeta_c-\zeta_b\in\R$
and
\begin{equation}\label{sim}
\overline{f(z)}=f(\overline z);
\end{equation}
thus,
$$
\overline{f(e^{i\theta})}=f(e^{i(2\pi-\theta)}), \ \ \theta\in[0,2\pi].
$$
Differentiating the latter formula with respect to $\theta$ and taking the modulus, yields 
\begin{equation}\label{f'}
|f'(e^{i\theta})|=|f'(e^{i(2\pi-\theta)})| , \ \ \theta\in[0,2\pi];
\end{equation}
thus, from (\ref{arclength}), we have that
$$
s(2\pi-\theta)=L-s(\theta),\hspace{1.5cm}\theta\in\R.
$$
From this formula and (\ref{f' g}), we obtain:
$$
\varphi(L-s(\theta))=\varphi(s(2\pi-\theta))
=\frac 1{2\pi|f'(e^{i(2\pi-\theta)})|},\hspace{1.5cm}\theta\in\R.
$$
Finally, from (\ref{f'}), it follows that
$$
\varphi(s)=\varphi(L-s),\hspace{1.5cm}s\in[0,L].
$$  

(ii) Suppose now $\varphi(s)=\varphi(L-s)$ for all $s\in\R.$ 
From (\ref{teta2}) we write
$$
\Phi(L-s)=2\pi\int_0^{L-s}\varphi(\sigma)d\sigma=
2\pi\int_s^L\varphi(L-\sigma)d\sigma=2\pi-\Phi(s),\ \ s\in[0,L].
$$
This property of $\Phi$ and (\ref{teta}) imply that
$$
\Phi(s(2\pi-\theta))=2\pi-\theta=2\pi-\Phi(s(\theta))=\Phi(L-s(\theta)),
$$
and hence
$$
s(2\pi-\theta)=L-s(\theta),\hspace{1.5cm}\theta\in\R,
$$
by the invertibility of $\Phi$. Then, by differentiating, we have that
$$
|f'(e^{i(2\pi-\theta)})|=s'(2\pi-\theta)=s'(\theta)=|f'(e^{i\theta})|
$$
for every $\theta\in\R$. Thus, by a change of variable and by simple
properties of the complex conjugate, we can write that, for $z\in D$,
\begin{eqnarray*}
\int_0^{2\pi}\frac{e^{it}+z}{e^{it}-z}\log|f'(e^{it})|dt
&=& \int_0^{2\pi}\frac{e^{it}+z}{e^{it}-z}\log|f'(e^{i(2\pi-t)})|dt \\
&=& \int_0^{2\pi}\frac{e^{i(2\pi-t)}+z}{e^{i(2\pi-t)}-z}\log|f'(e^{it})|dt \\
&=& \overline{\left(\int_0^{2\pi}\frac{e^{it}+\overline z}{e^{it}-\overline z}\log|f'(e^{it})|dt\right)}.
\end{eqnarray*}
Therefore, 
modulo a rotation, 
we have obtained that
\begin{equation*}\label{sim f'}
f'(z)=\overline{f'(\overline z)},\ \ z\in D,
\end{equation*}
and hence 
\begin{equation*}
f(z)=\overline{f(\overline z)},\ \ z\in D,
\end{equation*}
modulo a translation. Thus, $\mathcal{M}(\Omega)=\Omega$ for some reflection $\mathcal{M}.$
\end{proof}

%

\section{A formula involving curvature}

Recall that the curvature (with sign) $\kappa$ of a planar curve can be defined by the formula 
\begin{equation}\label{curvatura}
\kappa=\frac{d\psi}{ds},
\end{equation}
where $\psi$ is the angle between the positive real axis and the tangent (unit) vector.

By using the conformal map $f:D\rightarrow\Omega$ already introduced and
the Hilbert transform, we can express the curvature $\kappa$ of $\partial\Omega$ in terms
of the interior normal derivative $\varphi$ of the Green's function of $\Omega$. 

\begin{theorem} 
Let $\Omega\in\mathscr O$ and $\varphi$ be defined as usual.
Then $\varphi$ and the curvature $\kappa$ of $\partial\Omega$ are related by the formula:
\begin{equation}\label{curvaturaomega}
\kappa(s)=2\pi\varphi(s)\left[1-\frac 1{2\pi}\int_0^{|\partial\Omega|}
\cot\left(\frac{\Phi(s)-\Phi(\sigma)}2\right)\frac d{d\sigma}(\log\varphi)(\sigma)d\sigma\right],
\end{equation}
for $s\in[0,|\partial\Omega|]$, where $\Phi$ is defined as in (\ref{teta2}).
\end{theorem}

\begin{proof}
Let $f:D\rightarrow\Omega$ be as usual. Now we compute $\kappa$ in terms of $f.$
Define 
$$
\omega(\theta)=\arg(f'(e^{i\theta}))
$$
for $\theta\in[0,2\pi];$ the angle $\psi$ in (\ref{curvatura}) is given by
$$
\psi(\theta)=\arg\left(\frac d{d\theta}f(e^{i\theta})\right)=\omega(\theta)+\frac{\pi}2+\theta.
$$
From (\ref{curvatura}) and (\ref{fi}), we have that
\begin{equation}\label{curvaturafi}
\kappa(s)=\frac{d\psi}{d\theta}\frac{d\theta}{ds}=2\pi\varphi(s)[1+\omega'(\theta)],
\hspace{1.5cm}s\in[0,\partial\Omega].
\end{equation}
As is well-known (see \cite{Ko} and \cite{Pr}), since $\log|f'|$ and $\arg f'$ are the real and the imaginary
part of the analytic function $\log f',$ we have that 
\begin{equation}\label{argomento}
\arg f'(e^{i\theta})=\mathcal H(\log s')(\theta),
\end{equation}
being $s'(\theta)=|f'(e^{i\theta})|.$
Here, $\mathcal H$ is the 
Hilbert transformation on the unit circle, namely, 
$$
\mathcal H(\log s')(\theta)=\frac 1{2\pi}\int_0^{2\pi}\cot\left(\frac{\theta-t}2\right)\log(s'(t))dt.
$$
In our notations, (\ref{argomento}) can be rewritten as
$$
\omega=\mathcal H(\log s');
$$ 
thus,
$$
\omega'=\mathcal H(s''/s'),
$$ 
since $\mathcal H$ and $\displaystyle\frac{d}{d\theta}$ commute. From (\ref{curvaturafi}), we infer that
$$
\kappa(s(\theta))=2\pi\varphi\left[1+\mathcal H\left(\frac{s''}{s'}\right)(\theta)\right],
\ \ \theta\in[0,2\pi],
$$
and hence
$$
\kappa(s(\theta))=
2\pi\varphi(s(\theta))\left[1-\frac 1{2\pi}\int_0^{2\pi}
\cot\left(\frac{\theta-t}2\right)\frac{\varphi'(s(t))}{2\pi\varphi^2(s(t))}dt\right],
\ \ \theta\in[0,2\pi],
$$
from (\ref{fi}). Finally, we obtain (\ref{curvaturaomega}) by operating the change of variable
$\sigma=s(t)$ and by using (\ref{teta}).
\end{proof}

\begin{remark}
Let $\mathcal D2\mathcal N$ denote the Dirichlet-to-Neumann operator, that is 
$\mathcal D2\mathcal N$ maps the values on $\partial\Omega$ of any harmonic function in $\Omega$
to the values of its (interior) normal derivative on $\partial\Omega.$ Then, formula
(\ref{curvaturaomega}) can be rewritten as
$$
\kappa=2\pi\varphi[1+\mathcal D2\mathcal N(\log(\varphi))].
$$ 
\end{remark}

\section*{Acknowledgments} This research was partially supported by a PRIN
grant of the italian MIUR and by INdAM-GNAMPA.




\end{document}